\documentclass[a4paper,11pt]{article}
\usepackage{fullpage}
\usepackage{adjustbox}
\usepackage[rightcaption]{sidecap}
\usepackage{gensymb}

\usepackage[utf8]{inputenc}
\usepackage[english]{babel} 

\usepackage{amsfonts, amssymb, amsbsy, amsmath, amsthm}
\usepackage{dsfont}
\usepackage{authblk}
\allowdisplaybreaks

\usepackage{caption}
\usepackage[hyphens]{url}
\usepackage[hyphenbreaks]{breakurl}

\usepackage{mathtools}
\usepackage{xcolor}
\definecolor{darkblue}{rgb}{0,0,0.7}
\usepackage[colorlinks=true, allcolors=darkblue]{hyperref}
\usepackage[capitalise]{cleveref}
\Crefname{algocf}{Algorithm}{Algorithms}
\Crefname{equation}{Equation}{Equations}
\Crefname{figure}{Figure}{Figures}

\hypersetup{
    colorlinks=true,
    breaklinks=true,            %3
}

\usepackage{pgfplots}
\usepackage{caption}
\usepackage{wrapfig}
\usepackage{graphicx}
\usepackage{tikz-cd}
\usepackage{tikz,xcolor,xparse}
\usepackage{mathrsfs}
\usepackage{float}
\usepackage{animate}
\usetikzlibrary{shapes.geometric} 
\usetikzlibrary{positioning}
\usetikzlibrary{calc}
\usetikzlibrary{arrows.meta}
\usetikzlibrary{patterns}

\newtheorem{theorem}{Theorem}[section]

\usepackage{hhline}
\usepackage{stmaryrd}
\usepackage{xfrac}
\usepackage{algorithm}
\usepackage{algorithmic}
\numberwithin{equation}{section}

\newcommand{\EE}{\mathbb{E}}

\usepackage[textsize=tiny]{todonotes}

\begin{document}

\title{Extending the Continuum of Six-Colorings}
%Novel Six-Colorings for a Hadwiger-Nelson Variant
%
\author[1,2]{Konrad Mundinger}
\author[1,2]{Sebastian Pokutta}
\author[1,2]{Christoph Spiegel}
\author[1,2]{Max Zimmer}
\affil[1]{\small Technische Universit\"at Berlin, Institute of Mathematics}
\affil[2]{\small Zuse Institute Berlin, Department AIS2T, \emph{lastname}@zib.de}
\maketitle              % typeset the header of the contribution
\begin{abstract}
    We present two novel six-colorings of the Euclidean plane that avoid monochromatic pairs of points at unit distance in five colors and monochromatic pairs at another specified distance $d$ in the sixth color. Such colorings have previously been known to exist for $0.41 < \sqrt{2} - 1 \le d \le 1 / \sqrt{5} < 0.45$. Our results significantly expand that range to $0.354 \le d \le 0.657$, the first improvement in 30 years. Notably, the constructions underlying this were derived by formalizing colorings suggested by a custom machine learning approach.
\end{abstract}

\section{Introduction}\label{sec:introduction}

The Hadwiger-Nelson problem asks for the smallest number of colors needed to color the points of the Euclidean plane $\EE^2$ without any two points a unit distance apart having the same color. Viewing the plane as an infinite graph, with an edge between any two points if and only if the distance between them is $1$, motivates why this number is also referred to as the \emph{chromatic number of the plane} and denoted by $\chi (\EE^2)$. The problem goes back to 1950 and has since become one of the most enduring and famous open problems in combinatorial geometry and graph theory. For an extensive history of the problem and results related to it, we refer the reader to Jensen and Toft~\cite{jensen2011graph} as well as Soifer~\cite{nash2016open,soifer2009mathematical}.

By the de Bruijn–Erd\H{o}s theorem~\cite{bruijn1951colour}, and therefore assuming the axiom of choice, the problem is equivalent to finding the largest possible (vertex) chromatic number of a finite unit distance graph, that is a graph that can be embedded into the plane such that any two vertices are adjacent if and only if the corresponding points are at unit distance. Even without that connection, the chromatic number of any unit distance graph clearly gives a lower bound  for the chromatic number of the plane. The triangle is one obvious such graph, giving a lower bound of $3$, and the Moser spindle~\cite{moser1961solution} is the most famous example of a graph giving a lower bound of $4$. There had been no improvement to that lower bound since 1950 until de Grey famously established that $\chi (\EE^2) \ge 5$ through a graph of order $1581$ in 2018~\cite{DeGrey2018ChromaticNumber}. Simplifying and reducing the size of this construction has been of great interest to the extent of being the topic of a Polymath project~\cite{Heule2018ComputingSmallGraphs,exoo2020chromatic, parts2020graph, Polymath2021ChromaticNumber}.

Regarding upper bounds, there is a large number of distinct $7$-colorings of the plane that avoid monochromatic pairs at unit distance, the first of which (using a tiling of the plane with congruent regular hexagons) was already observed back in 1950 by Isbell~\cite{nash2016open,soifer2009mathematical}. This upper bound of $\chi (\EE^2) \le 7$ has remained unchanged since and many variants of the original question have therefore been proposed in the hopes of shedding some light on why this problem has proven so stubborn. To state one such variant, we say that an $n$-coloring of the plane has \emph{coloring type} $(d_1, \ldots, d_n)$ if color $i$ does not realize distance $d_i$~\cite{soifer1992relatives,soifer1992six}. This gives a measurement of how close this coloring is to achieving the original goal and can be seen as a defining a natural {\lq}off-diagonal{\rq} variant of the original problem. Finding a coloring of type $(1,1,1,1,1,1)$ would obviously improve the upper bound of $\chi (\EE^2)$ to $6$.

Stechkin found a coloring of type $(1,1,1,1,1/2,1/2)$, which was published by Raiskii in 1970~\cite{raiskii1970realization}, and Woodall found a coloring of type $(1, 1, 1, 1/\sqrt{3}, 1/ \sqrt{3}, 1/ \sqrt{12})$ in 1973~\cite{woodall1973distances}. The first six-coloring to feature a non-unit distance in only one color has type $(1,1,1,1,1,1/\sqrt{5})$ and was found by Soifer in 1991~\cite{soifer1992six}. Hoffman and Soifer also found a coloring of type $(1,1,1,1,1,\sqrt{2} - 1)$ in 1993~\cite{hoffman1993almost, hoffman1996another}. Both of these constructions are in fact part of a family that realizes  $(1,1,1,1,1,d)$ for any $ \sqrt{2} - 1 \le d \le 1 / \sqrt{5}$~\cite{hoffman1996another, soifer1994infinite, soifer2009mathematical}, leading Soifer~\cite{Soifer1994SixRealizable} to pose the “still open and extremely difficult”~\cite{nash2016open} problem of determining the \emph{continuum of six colorings} $X_6$, that is the set of all $d$ for which there exists a six-coloring of the plane of type $(1,1,1,1,1,d)$. To the best of our knowledge, no improvements have been suggested in the last 30 years.

We propose two novel six-colorings of the plane, one parameterized by $d$ and the other valid for a range of values for $d$ without modifications, that together significantly expand the range of $d$ known to be in $X_6$. The first is a valid coloring of type $(1,1,1,1,1,d)$ as long as $0.354 \le d \le 0.553$ and the second covers the range of $0.418 \le d \leq 0.657$.% The fact that the second construction is \emph{not} parameterized by $d$ also makes it a valid coloring for some specific periodicic intervals when $d > 1$.

\begin{theorem}\label{thm:main}
    $X_6$ contains the closed interval $[0.354, 0.657]$.
\end{theorem}

It should be noted that both constructions were derived by formalizing colorings that were suggested by a custom machine learning approach in which a Neural Network was trained to represent a coloring of a specified type or range of types. We will briefly touch upon this in \cref{sec:discussion} and otherwise go into more detail about this approach and potential other applications in a separate publication. This work is intended to give a formal justification of \cref{thm:main}, with the first coloring being explored in \cref{sec:firstcoloring} and the second in \cref{sec:secondcoloring}.

\begin{figure}[htbp]
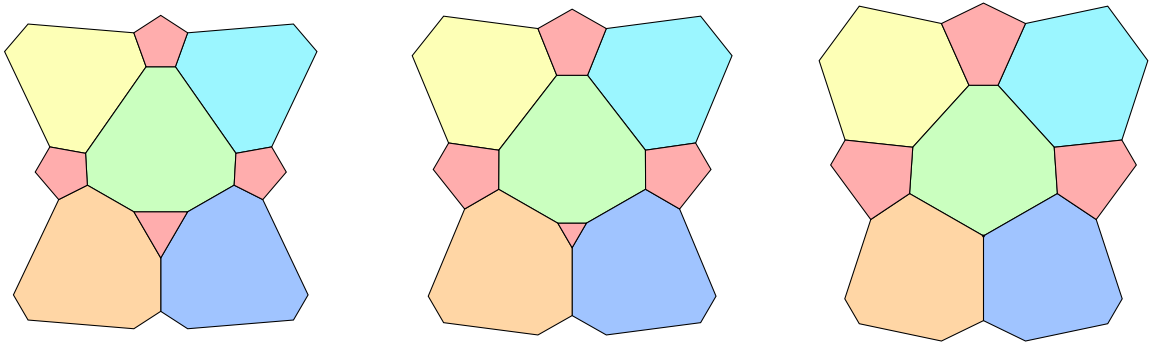

    \centering
    % First Block
    \begin{minipage}{.32\textwidth}
        \centering
        \pgfmathsetmacro{\avoid}{0.354}
        \input{figures/firstcoloring/block} % Adjust path as necessary
    \end{minipage}\hfill % The \hfill command adds horizontal space between the minipages
    % Second Block
    \begin{minipage}{.32\textwidth}
        \centering
        \pgfmathsetmacro{\avoid}{0.45}
        \input{figures/firstcoloring/block} % Adjust path as necessary
    \end{minipage}\hfill
    % Third Block
    \begin{minipage}{.32\textwidth}
        \centering
        \pgfmathsetmacro{\avoid}{0.553}
        \input{figures/firstcoloring/block} % Adjust path as necessary
    \end{minipage}
    \caption{Building block of first coloring for $d = 0.354$, $0.45$, and $0.553$.} \label{fig:firstcoloring_block}
\end{figure}

\section{A construction for $0.354 \le d \le 0.553$}\label{sec:firstcoloring}

The first constructions is made up of four different polytopal shapes, a detailed description of which is given in the appendix. The equidiagonal pentagon and the equilateral triangle respectively described in \cref{fig:first_pentagon} and \cref{fig:first_triangle} together are colored with the sixth color (red) in which we are avoiding points at distance $d$. The octagons described in \cref{fig:first_octagon} receive three of the other five colors (orange, green, and blue) and the hexagons described in \cref{fig:first_hexagon} receive the remaining two (yellow and turquoise). All shapes are uniquely parameterized by the choice of $d$, with the exception of the pentagon, which has an additional degree of freedom in the form of $\alpha_1$, the angle at the {\lq}top{\rq} of the pentagon. We will later determine the range of valid $\alpha_1$ depending on $d$ numerically and see that this additional variable can be fixed by linearly interpolating between two extremal values (though other options can also be valid depending on $d$).

\begin{figure}[htbp]
    \centering
    \input{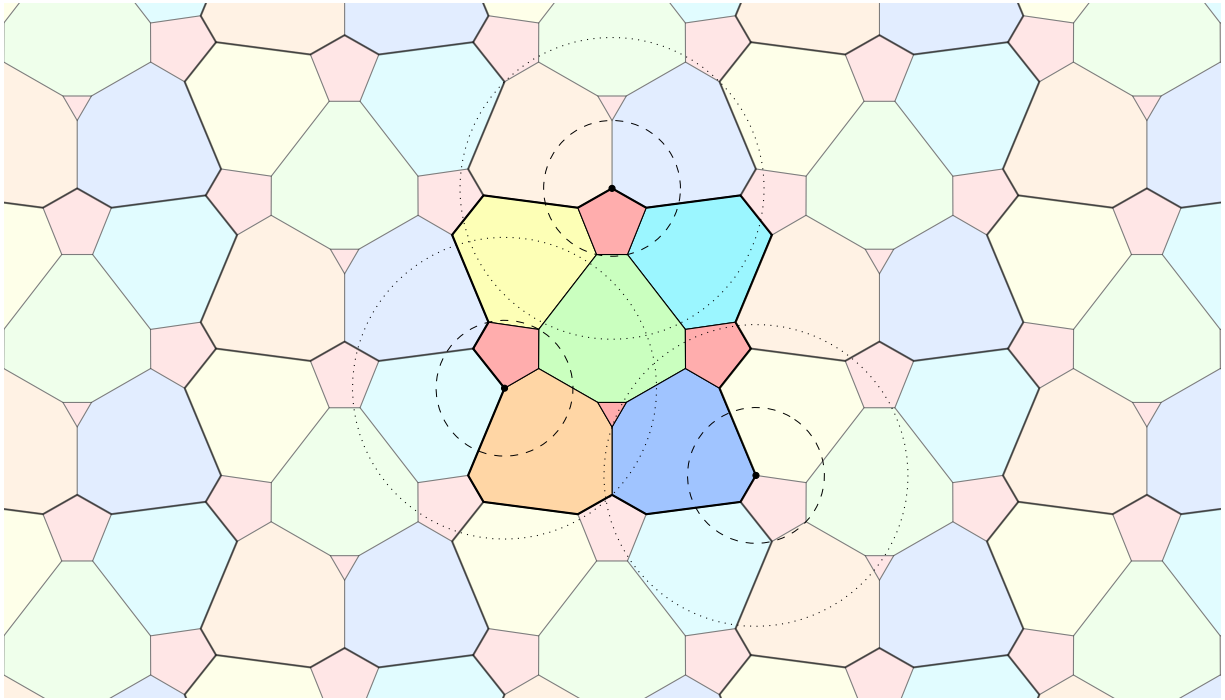}
    \caption{Illustration of the first coloring with circles at unit distance (dotted) and distance $d$ (dashed) highlighted at three critical points. %An animated version can be found in \cref{fig:firstcoloringanimated} in the appendix.
    }
    \label{fig:firstcoloring}
\end{figure}

A copy of three pentagons, one triangle, three octagons and two hexagons together form the building block of the first coloring that is illustrated in \cref{fig:firstcoloring_block} for three different valid choices of $d$. Note that the triangle disappears as $d$ approaches the upper end of the valid spectrum. Looking at the overall construction in \cref{fig:firstcoloring}, it is visually clear that the only conditions that are at risk making this construction invalid are given be the following set of constraints:

\vspace{-1em}
\begin{minipage}[t]{0.4\textwidth}
\begin{align}
    s_4 & \le d \label{eq:cond1} \\
    s_5 & \ge d \label{eq:cond2} \\
    w_1 & \le 1 \label{eq:cond3}
\end{align}
\end{minipage}%
\begin{minipage}[t]{0.5\textwidth}
\begin{align}
    w_2 & \le 1 \label{eq:cond4} \\
    w_3 & \le 1 \label{eq:cond5} \\ 
    h_1 + h_3 + d & \ge 1 \label{eq:cond6}
\end{align}
\end{minipage}%
\vspace{1.5em}

Here $h_1$ is the height of the pentagon, $h_3$ the height of the triangle, $s_4$ is the side length of the triangle, $s_5$ the length of the longest side of the octagon, $w_1$ and $w_2$ two different widths of octagon, and $w_3$ the width of the hexagon. Note that a more detailed description of variables alongside the corresponding shape is given in the appendix. \cref{fig:firstcolorincritical} in the appendix also gives a visual representation of the six cases.

Unfortunately we were unable to derive a closed form expression for the range of $d$ for which a valid choice of $\alpha_1$ can be found. However, it is easy to numerically verify that for $d \in [0.354, 0.553]$ such a choice can be made. Furthermore, by linearly interpolating between the two extreme points, that is by choosing $\alpha_1 = 113.7 + (d - 0.354) \, 14.11 / 0.299$, we can remove the additional degree of freedom in the definition of the pentagon. Finally, we note that there is always an appropriate choice for the color on the boundaries between the shapes.

\begin{figure}[htbp]
    \centering
    \input{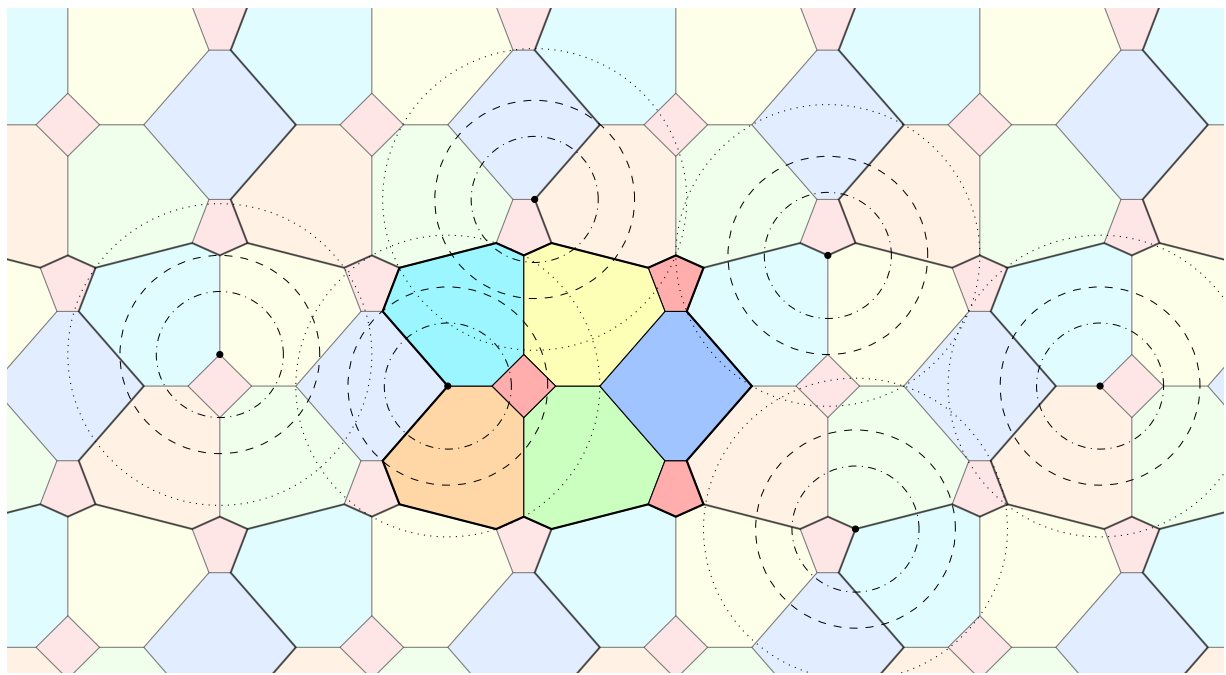}
    \caption{Illustration of the second coloring with circles at unit distance (dotted), and distance $d_{\max}$ (dashed), and distance distance $d_{\min}$ (dash-dotted) highlighted at six critical points.}
    \label{fig:secondcoloring}
\end{figure}

\section{A construction for $0.418 \le d \le 0.657$}\label{sec:secondcoloring}

Let $d_{\max}$ be the real root of $d^4 + 5 \sqrt{3} d^3 + 18 d^2 - 3 \sqrt{3} d - 7 = 0$ closest to $0.65$ and $d_{\min} = \sqrt{3} - 2 \, d_{\max}$. Note that a closed form for $d_{\max}$ is given by

\begin{align*}
    d_{\max} & = -(5 \sqrt{3})/4 + 1/2 \, \big(27/4 + 1/3 \, (7290 - 15 \sqrt{1821})^{1/3} + (5 \, (486 + \sqrt{1821}))^{1/3}/3^{2/3} \big)^{1/2} \\
    & \quad + 1/2 \, \Big(27/2 - 1/3 \, (7290 - 15 \sqrt{1821})^{1/3} - (5 (486 + \sqrt{1821}))^{1/3}/3^{2/3} \\
    & \quad + 9/4 \, \big( 3/(27/4 + 1/3 \, (7290 - 15 \sqrt{1821})^{1/3} + (5 (486 + \sqrt{1821}))^{1/3}/3^{2/3}))\big)^{1/2} \Big)^{1/2}.
\end{align*}
We can easily verify numerically that $d_{\min} \le 0.418 \le d \le 0.657 \le d_{\max}$ and the second construction  will in fact be valid for any $d \in [d_{\min}, d_{\max}]$. It is again made up of four different polytopal shapes, a detailed description of which is given in the appendix. The pentagon and square described in \cref{fig:second_pentagon} together are colored with the sixth color (red) in which we are avoiding points at distance $d$. The heptagon described in \cref{fig:second_heptagon} receives four of the other five colors (orange, green, yellow, and turquoise) while hexagon described in \cref{fig:second_heptagon} receives the last remaining color (blue). A copy of two pentagons, one square, four heptagons and one hexagon together form the building block of the second coloring, which is illustrated in \cref{fig:secondcoloring}.

% \textbf{Discuss for which intervals this holds}

% \begin{figure}[htbp]
%     \centering
%     \input{figures/secondcoloring/completelarge}
%     \caption{....}
%     \label{fig:secondcoloringlarge}
% \end{figure}

\section{Discussion and Outlook}\label{sec:discussion}

We conclude by noting that there was a significant technical component to these new constructions. We developed a custom machine learning approach in which we had a Neural Network represent a (probabilistic) six-coloring of the plane. We defined a loss function based on the likelihood that two points at unit distance (or at distance $d$) are monochromatic with respect to the right color(s) and updated the parameters of the neural network to minimize this loss using a variant of stochastic gradient descent. The resulting output was detailed enough to inspire the above constructions, though formally describing them and verifying their correctness still required a fair amount of manual effort.

This is significant, as the opposite problem, that is finding (unit) distance graphs with large chromatic number, has recently been heavily dominated by computer-based (SAT) approaches. As far as we are aware, the only previous documented attempt of using computers to derive colorings was limited to using SAT solvers to essentially obtain very low-resolution images in the context of the Polymath project. Of course, we also attempted to derive a six-coloring when $d = 1$, but were unable to find anything beyond the almost-successful six-coloring of Pritikin~\cite{pritikin1998all}, which is closely connected to a coloring of Pegg~\cite{soifer2009mathematical} and was recently improved by Parts~\cite{parts2020percent}. We believe that this adds some more weight to the commonly held belief that $\chi(\mathbb{E}^2) = 7$.  We were also unable to find a five-coloring of type $(1,1,1,1,d)$, indicating that the \emph{almost chromatic number}~\cite{soifer2009mathematical} perhaps satisfies $\chi_a(\mathbb{E}^2) = 6$.

We conclude by noting that we do not believe our constructions to be fully optimized in terms of which values of $d$ in $X_6$ are covered. As the work of Parts~\cite{parts2020percent} on \emph{almost} coloring the plane with six colors while avoiding points at unit distance showed, a description of these colorings can become incredibly nuanced and intricate. Connected to this, we note that there were signs in our computational results that, with some more technical adjustments, the second coloring could also be valid for values $d > 1$ or even periodically for certain intervals.

\paragraph{Acknowledgements} Research reported in this paper was partially supported through the German Federal Ministry of Education and Research (fund number 01IS23025B) as well as the Deutsche Forschungsgemeinschaft (DFG, German Research Foundation) under Germany’s Excellence Strategy – The Berlin Mathematics Research Center MATH+ (EXC-2046/1, project ID: 390685689).

% \bibliographystyle{splncs04}
% \bibliography{bib}
%

\newpage

\appendix

\section{Critical cases in the first coloring}

\begin{figure}[H]
    \centering
    \input{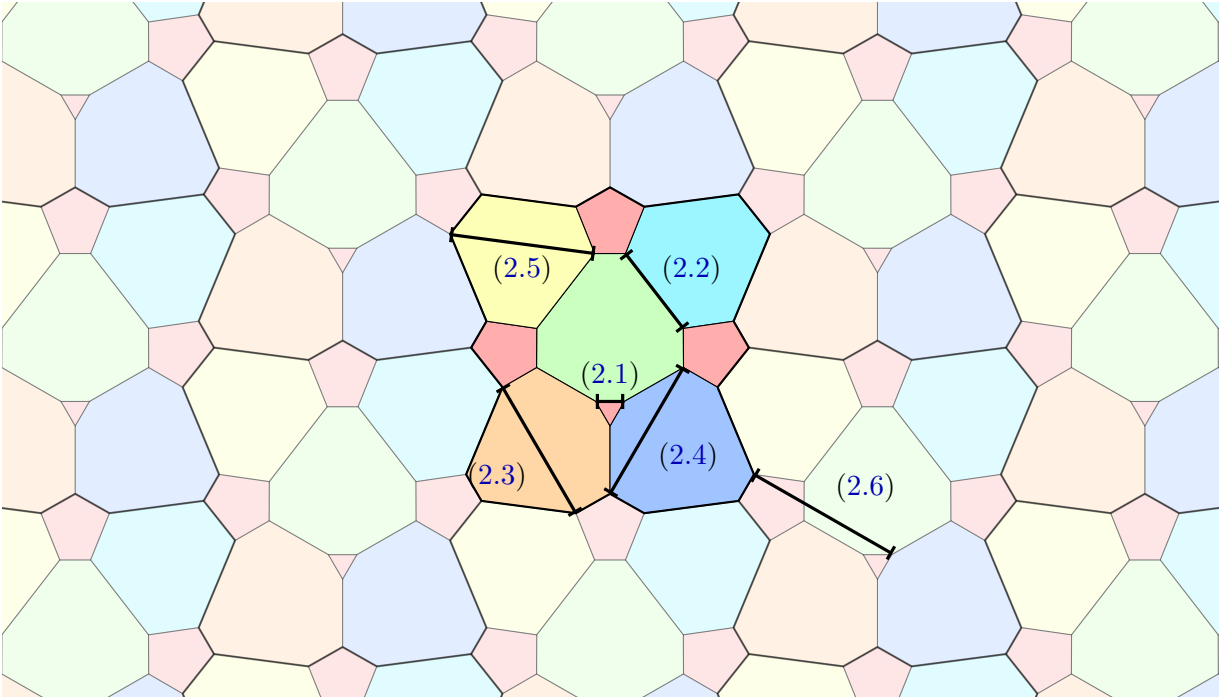}
    \caption{Illustrating the critical cases \eqref{eq:cond1} to \eqref{eq:cond6} in the first coloring.}
    \label{fig:firstcolorincritical}
\end{figure}

% \section{Animation of the first coloring}

% \begin{figure}[H]
%     \centering
%     \input{figures/firstcoloring/animation}
%     \caption{Animated illustration of the first coloring with circles at unit distance (dotted) and distance $d$ (dashed) highlighted at three critical points. Note that only some PDF viewers support displaying TikZ animations.}
%     \label{fig:firstcoloringanimated}
% \end{figure}

\section{Building blocks of the first coloring}

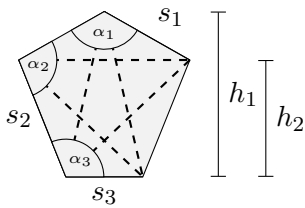
\begin{figure}[H]
    \begin{minipage}[t]{0.5\textwidth}
    \centering
    \vspace{4em}
    \tikzset{every picture/.style={scale=5}}

\pgfmathsetmacro{\avoid}{0.45}
\input{figures/firstcoloring/variables}

\begin{tikzpicture}
    \coordinate (A) at (0,0);
    \coordinate (B) at ([shift=(180+90-\pat/2:\pst)]A);
    \coordinate (C) at ([shift=(180+90+\pat/2:\pst)]A);
    
    \coordinate (D) at ([shift=(180+90-\pait/2:\avoid)]A);
    \coordinate (E) at ([shift=(180+90+\pait/2:\avoid)]A);

    \draw[fill=black!5] (A) -- (B) -- (D) -- (E) -- (C) -- cycle;
    
    \draw[thick, dashed] (B) -- (C);
    \draw[thick, dashed] (A) -- (E);
    \draw[thick, dashed] (A) -- (D);
    \draw[thick, dashed] (C) -- (D);
    \draw[thick, dashed] (B) -- (E);
        
    \draw [{Bar[width=2mm]}-{Bar[width=2mm]}] ($(A) + (0.3,0)$) -- ++(0,-\ph) node[pos=0.5, right] {$h_1$};
    
    \draw [{Bar[width=2mm]}-{Bar[width=2mm]}] ($(C) + (0.2,0)$) -- ++(0,{-\phb}) node[pos=0.5, right] {$h_2$};

    \draw  (A) -- (C) node[pos=0.5, above right] {$s_1$};
    \draw  (B) -- (D) node[pos=0.5, left] {$s_2$};
    \draw  (D) -- (E) node[pos=0.5, below] {$s_3$};
    
    \draw[fill=black!5] (A) -- ++(270+\pat/2:0.1) arc (270+\pat/2:270-\pat/2:0.1) node[pos=0.5, above] {\tiny  $\alpha_1$}  -- cycle;
    
    \draw[fill=black!5] (B) -- ++(\past:0.1) arc (\past:\past-\pas:0.1) node[pos=0.4, left, inner sep=2pt] {\tiny $\alpha_2$}  -- cycle;
    
    \draw[fill=black!5] (D) -- ++(0:0.1) arc (0:0+\pab:0.1) node[pos=0.6, below, inner sep=5pt] {\tiny $\alpha_3$}  -- cycle;
    
    % \draw[fill=black] (A) circle (0.01);
    % \draw[line width=0.1pt] (A) ++(250:\avoid) arc (250:290:\avoid);

    % \draw[fill=black] (B) circle (0.01);
    % \draw[line width=0.1pt] (B) ++(310:\avoid) arc (310:370:\avoid);
    
    % \draw[fill=black] (E) circle (0.01);
    % \draw[line width=0.1pt] (E) ++(90:\avoid) arc (90:145:\avoid);
\end{tikzpicture}
    \end{minipage}%
    \begin{minipage}[t]{0.5\textwidth}
    \begin{align*}
        s_1 & = d/2 \, \csc(\alpha_1 / 2) \\
        t_1 & = 2 \arccos \big(\csc(\alpha_1 / 2)/4 \big) - \alpha_1 \\
        s_3 & = 2d \, \sin(t_1/2) \\
        h_1 & = d \, \cos(t_1/2) \\ 
        h_2 & = h_1 - (d/2) \, \cot(\alpha_1/2) \\
        s_2 & = \sqrt{ h_2^2 + (d-s_3)^2 / 4} \\
        \alpha_2 & = 90\degree - \alpha_1/2 + \arcsin(h_2 / s_2)\\
        \alpha_3 & = 270\degree - \alpha_1/2 - \alpha_2
    \end{align*}
    \end{minipage}
    \caption{An equidiagonal pentagon with each diagonal of length $d$, highlighted by dashed lines, used for the red color avoiding points at distance $d$ in the first coloring.
    }
    \label{fig:first_pentagon}
\end{figure}

\begin{figure}[H]
    \begin{minipage}[t]{0.35\textwidth}
    \centering
    \vspace{3em}
    \tikzset{every picture/.style={scale=5}}                           

\pgfmathsetmacro{\avoid}{0.45}
\input{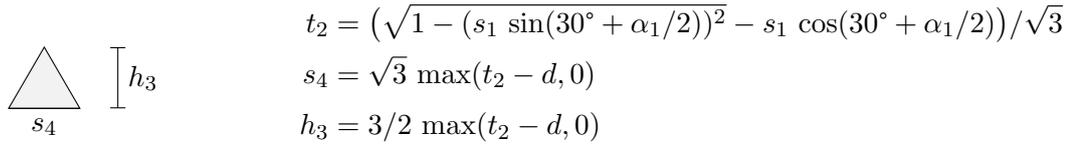}

\begin{tikzpicture}
    \coordinate (A) at ([shift=(90:\tc)](0,0);
    \coordinate (B) at ([shift=(-30:\tc)](0,0);
    \coordinate (C) at ([shift=(-150:\tc)](0,0);
    
    \draw[fill=black!5] (A) -- (B) -- (C) -- cycle; 
    
    \draw  (B) -- (C) node[pos=0.5, below] {$s_4$};
    
    \draw [{Bar[width=2mm]}-{Bar[width=2mm]}] ($(B)
     + (0.1,0)$) -- ++(0,\th) node[pos=0.5, right] {$h_3$};
\end{tikzpicture}
    \end{minipage}%
    \begin{minipage}[t]{0.6\textwidth}
    \begin{align*}
        t_2 & = \big( \sqrt{1 - (s_1 \, \sin(30\degree + \alpha_1/2))^2} - s_1 \, \cos(30\degree+\alpha_1/2) \big) / \sqrt{3} \\
        s_4 & = \sqrt{3} \, \max(t_2-d,0) \\
        h_3 & = 3/2 \, \max(t_2-d,0)
    \end{align*}
    \end{minipage}
    \caption{An equilateral triangle, used for the red color avoiding points at distance $d$ in the first coloring.
    }
    \label{fig:first_triangle}
\end{figure}

\begin{figure}[H]
    \begin{minipage}[t]{0.5\textwidth}
    \centering
    \vspace{1em}
    \tikzset{every picture/.style={scale=3}}

% Distance avoided in last color
\pgfmathsetmacro{\avoid}{0.45}                                
\input{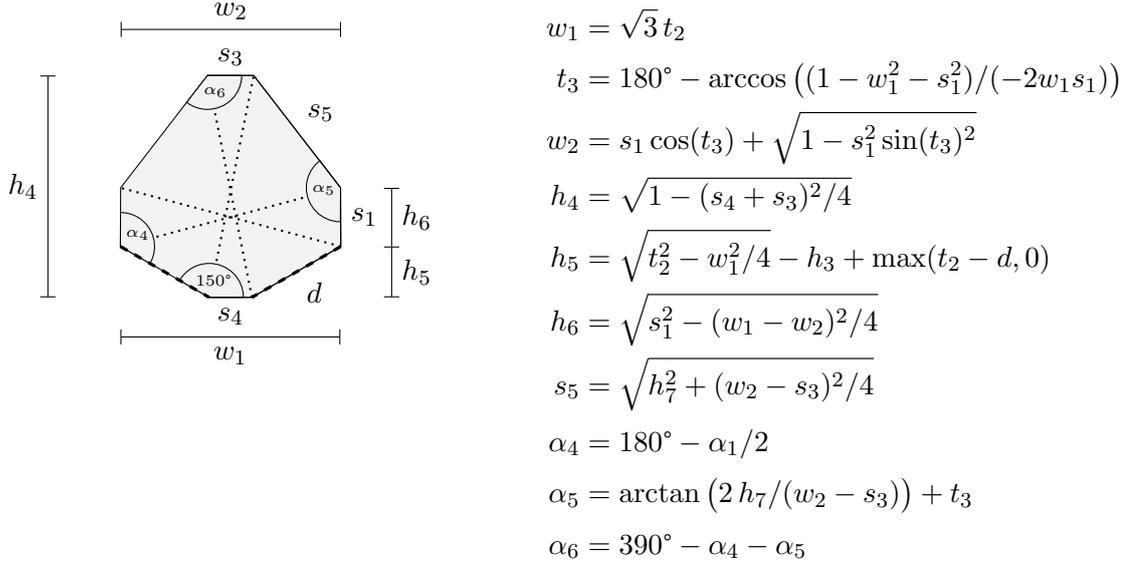}

\begin{tikzpicture}
    \coordinate (A) at (0,0);
    
    \coordinate (B) at ([shift=(150:\tc)]A);
    \coordinate (C) at ([shift=(30:\tc)]A);
    
    \coordinate (BB) at ([shift=(150:\bd)]A);
    \coordinate (CC) at ([shift=(30:\bd)]A);
    
    \coordinate (BBB) at ([shift=(150-\pat/2:\pst)]BB);
    \coordinate (CCC) at ([shift=(30+\pat/2:\pst)]CC);

    \coordinate (E) at ([shift=(\bait:1)]B);
    \coordinate (D) at ([shift=(180-\bait:1)]C);
    
    \draw[fill=black!5] (E) -- (CCC) -- (CC) -- (C) -- (B) -- (BB) -- (BBB) -- (D) -- cycle;

    \draw[thick, dotted] (BB) -- (CCC);
    \draw[thick, dotted] (CC) -- (BBB);
    
    \draw[thick, dotted] (B) -- (E);
    \draw[thick, dotted] (C) -- (D);
    
    \draw[fill=black!5] (B) -- ++(0:0.15) arc (0:\bat:0.15) node[pos=0.55, below, inner sep=4pt] {\tiny $150\degree$}  -- cycle;

    \draw[fill=black!5] (BB) -- ++(-30:0.15) arc (-30:-30+\bast:0.15) node[pos=0.5, below left, inner sep=0pt] {\tiny $\alpha_4$}  -- cycle;

    \draw[fill=black!5] (CCC) -- ++(-30-\pat/2:0.15) arc (-30-\pat/2:-30-\pat/2-\basb:0.15) node[pos=0.6, right, inner sep=2pt] {\tiny $\alpha_5$}  -- cycle;
    
    \draw[fill=black!5] (D) -- ++(0:0.15) arc (0:-\bab:0.15) node[pos=0.6, above, inner sep=4pt] {\tiny $\alpha_6$}  -- cycle;
    
    % \draw[fill=black] (CC) circle (0.01);
    % \draw[thick] (CC) ++(160:1) arc (160:180:1);
    % \draw[line width=0.1pt] (CC) ++(195:\avoid) arc (195:220:\avoid);
    
    % \draw[fill=black] (CCC) circle (0.01);
    % \draw[thick] (CCC) ++(-155:1) arc (-155:-180:1);

    % \draw[fill=black] (C) circle (0.01);
    % \draw[thick] (C) ++(90:1) arc (90:110:1);
    % \draw[line width=0.1pt] (C) ++(20:\avoid) arc (20:40:\avoid);

    \draw[very thick, dashed] (BB) -- (B);
    \draw[very thick, dashed] (CC) -- (C);

    \draw  (B) -- (C) node[pos=0.5, below] {$s_4$};
    \draw  (E) -- (D) node[pos=0.5, above] {$s_3$};
    \draw  (CC) -- (CCC) node[pos=0.5, right] {$s_1$};
    \draw  (E) -- (CCC) node[pos=0.5, above right] {$s_5$};
    \draw  (C) -- (CC) node[pos=0.5, below right] {$d$};
    
    \draw [{Bar[width=2mm]}-{Bar[width=2mm]}] ($(D) + (-0.7,0)$) -- ++(0,-\bh) node[pos=0.5, left] {$h_4$};

    % \draw [{}-{Bar[width=2mm]}] ($(B) + (0.8,\bht+\bhm)$) -- ++(0,\bhb) node[pos=0.5, right] {$h_7$};
    \draw [{}-{Bar[width=2mm]}] ($(B) + (0.8,\bht)$) -- ++(0,\bhm) node[pos=0.5, right] {$h_6$};
    \draw [{Bar[width=2mm]}-{Bar[width=2mm]}] ($(B) + (0.8,0)$) -- ++(0,\bht) node[pos=0.5, right] {$h_5$};

    \draw [{Bar[width=2mm]}-{Bar[width=2mm]}] ($(BB) + (0,-0.4)$) -- ++(\bwb,0) node[pos=0.5, below] {$w_1$};

    \draw [{Bar[width=2mm]}-{Bar[width=2mm]}] ($(BBB) + (0,0.7)$) -- ++(\bwt,0) node[pos=0.5, above] {$w_2$};    
    
\end{tikzpicture}
    \label{fig:octagon}
    \end{minipage}%
    \begin{minipage}[t]{0.5\textwidth}
    \begin{align*}
        w_1 & = \sqrt{3} \, t_2 \\
        t_3 & = 180\degree - \arccos \big( (1-w_1^2-s_1^2)/(-2w_1 s_1) \big) \\
        w_2 & = s_1 \cos(t_3)+\sqrt{1-s_1^2 \sin(t_3)^2} \\
        h_4 & = \sqrt{1 - (s_4+s_3)^2 / 4} \\
        h_5 & = \sqrt{t_2^2-w_1^2/4} - h_3 + \max(t_2-d,0)\\
        h_6 & = \sqrt{s_1^2 - (w_1-w_2)^2 / 4}\\ 
        % h_7 & = h_4 - h_5 - h_6 \\
        s_5 & = \sqrt{ h_7^2 + (w_2-s_3)^2 / 4} \\
        \alpha_4 & = 180\degree - \alpha_1/2 \\
        \alpha_5 & =  \arctan\big(2 \, h_7 / (w_2-s_3) \big)+t_3\\
        \alpha_6 & =  390\degree-\alpha_4-\alpha_5 
    \end{align*}
    \end{minipage}
    \caption{An axisymmetric octagon in which four of the diagonals have unit length, highlighted by dotted lines, and two of the sides have length $d$, highlighted by dashed lines. Used for the orange, green and blue color avoiding points at unit distance in the first coloring.}\label{fig:first_octagon}
\end{figure}

\begin{figure}[H]
    \begin{minipage}[t]{0.5\textwidth}
    \centering
    \vspace{1em}
    \tikzset{every picture/.style={scale=3}}

% Distance avoided in last color
\pgfmathsetmacro{\avoid}{0.45}
\input{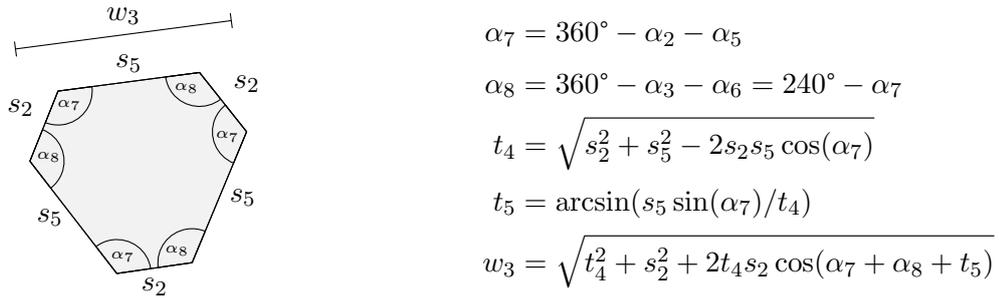}

\begin{tikzpicture}
        \coordinate (A) at (0,0);
        \coordinate (B) at ([shift=(180-\pab:\pss)]A);
        \coordinate (C) at ([shift=(-\pab+\faf:\bsb)]B);
        \coordinate (D) at ([shift=(-\pab+\faf-180+\fas:\pss)]C);
        \coordinate (E) at ([shift=(-\pab+2*\faf+\fas:\bsb)]D);
        \coordinate (F) at ([shift=(-\pab+2*\faf+2*\fas-180:\pss)]E);

        \draw[fill=black!5]  (A) -- (B) -- (C) -- (D) -- (E) -- (F) --cycle;

        \draw  (A) -- (B) node[pos=0.5, above left] {$s_2$};
        \draw  (B) -- (C) node[pos=0.5, above] {$s_5$};
        \draw  (C) -- (D) node[pos=0.5, above right] {$s_2$};
        \draw  (D) -- (E) node[pos=0.5, right] {$s_5$};
        \draw  (E) -- (F) node[pos=0.5, below] {$s_2$};
        \draw  (F) -- (A) node[pos=0.5, left] {$s_5$};
        
        % \draw[fill=black] (A) circle (0.01);
        % \draw[thick] (A) ++(0:1) arc (0:20:1);
        
        % \draw[fill=black] (B) circle (0.01);
        % \draw[thick] (B) ++(-90:1) arc (-90:10:1);
        
        % \draw[fill=black] (C) circle (0.01);
        % \draw[thick] (C) ++(-130:1) arc (-130:-90:1);
        
        % \draw[fill=black] (D) circle (0.01);
        % \draw[thick] (D) ++(-200:1) arc (-200:-160:1);
        
        % \draw[fill=black] (E) circle (0.01);
        % \draw[thick] (E) ++(100:1) arc (100:150:1);
        
        % \draw[fill=black] (F) circle (0.01);
        % \draw[thick] (F) ++(60:1) arc (60:100:1);
            
        % \draw (A) -- ([shift=(180-\pab-\fad:\fd)]A);
        % \draw[ultra thin, color=blue] (A) -- ++(180-\pab:0.1) arc (180-\pab:180-\pab-\fad:0.1) -- cycle;
        
        % \draw[ultra thin, color=blue] (A) -- ++(180-\pab:0.11) arc (180-\pab:180-\pab-\fas:0.11) -- cycle;
                
        % \draw[ultra thin, color=blue] (B) -- ++(-\pab:0.1) arc (-\pab:-\pab+\faf:0.1) -- cycle;
        
        % \draw[ultra thin, color=blue] (C) -- ++(360-\pab-\fad:0.1) arc (360-\pab-\fad:360-\pab-\fad-\fadd:0.1) -- cycle;
        
        % \draw[ultra thin, color=blue] (C) -- ++(360-\pab-\fad:0.11) arc (360-\pab-\fad:360-\pab-\fad+\faddd:0.11) -- cycle;

        % \draw[dotted] (A) circle (\fw);

        \coordinate (X) at ([shift=(90-\pab+\faf:0.5)]A);
        \draw [{Bar[width=2mm]}-{Bar[width=2mm]}] (X) -- ++(-\pab+\faf:\fw) node[pos=0.5, above] {$w_3$};

        % \draw (A) circle (\fw);
        % \draw (B) circle (\fw);
        % \draw (C) circle (\fw);

        % \coordinate (X) at ([shift=(180-\pab+\faf:0.5)]B);
        % \draw [{Bar[width=2mm]}-{Bar[width=2mm]}] (B) -- ++(-\pab+\faff:\fw) node[pos=0.5, above] {$w_4$};

        % faf
        \draw (D) -- ++(-\pab+2*\faf+\fas:0.15) arc (-\pab+2*\faf+\fas:-\pab+\faf+\fas:0.15) node[pos=0.5, inner sep=1pt, right] {\tiny $\alpha_7$}  -- cycle;
        
        \draw (B) -- ++(-\pab+\faf:0.15) arc (-\pab+\faf:-\pab:0.15) node[pos=0.4, inner sep=1pt, above left] {\tiny $\alpha_7$}  -- cycle;

        \draw (F) -- ++(-\pab+2*\faf+2*\fas-180:-0.15) arc (-\pab+2*\faf+2*\fas-180:-\pab+2*\faf+2*\fas-180+\faf:-0.15) node[pos=0.6, below, inner sep=4pt] {\tiny $\alpha_7$}  -- cycle;

        % fas
        \draw (E) -- ++(-\pab+2*\faf+2*\fas-180:0.15) arc (-\pab+2*\faf+2*\fas-180:-\pab+2*\faf+2*\fas-180-\fas:0.15) node[pos=0.6, below, inner sep=5pt] {\tiny $\alpha_8$}  -- cycle;
    
        \draw (C) -- ++(-\pab+\faf-180+\fas:0.15) arc (-\pab+\faf-180+\fas:-\pab+\faf-180+\fas-\fas:0.15) node[pos=0.5, inner sep=4pt, above] {\tiny $\alpha_8$}  -- cycle;

        \draw (A) -- ++(180-\pab:0.15) arc (180-\pab:180-\pab-\fas:0.15) node[pos=0.5, left, inner sep=1pt] {\tiny $\alpha_8$}  -- cycle;
    
\end{tikzpicture}
    \end{minipage}%
    \begin{minipage}[t]{0.5\textwidth}
    \begin{align*}
        \alpha_7 & = 360\degree - \alpha_2 - \alpha_5 \\
        \alpha_8 & = 360\degree - \alpha_3 - \alpha_6  = 240\degree - \alpha_7 \\
        t_4 & = \sqrt{s_2^2 + s_5^2 - 2 s_2 s_5 \cos(\alpha_7)} \\
        t_5 & = \arcsin(s_5 \sin(\alpha_7)/t_4) \\
        w_3 & = \sqrt{t_4^2 + s_2^2 + 2 t_4 s_2 \cos(\alpha_7+\alpha_8+t_5)}
    \end{align*}
    \end{minipage}
    \caption{A hexagon with two angles and two side lengths. Used for the yellow and turqouise color avoiding points at unit distance in the first coloring. Note that it is in general \emph{not} axisymmetric.}
    \label{fig:hexagon}\label{fig:first_hexagon}
\end{figure}

\section{Building blocks of the second coloring}

\begin{figure}[H]
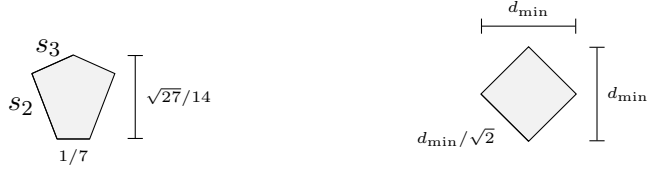

    \centering
    \begin{minipage}[t]{0.35\textwidth}
    \centering
    \vspace{3.0em}
    \pgfmathsetmacro{\scale}{3}
\tikzset{every picture/.style={scale=\scale}}

\input{figures/secondcoloring/variables}

\begin{tikzpicture}
        \coordinate (PA) at (1/2, {sqrt(3)/2});
        \coordinate (I2) at ({-1/4 + \auxw * cos(\auxa)}, {sqrt(3)/4 + \maxavoid/2 + \auxw * sin(\auxa)});
        \coordinate (I3) at ({5/4 - \auxw * cos(\auxa)}, {sqrt(3)/4 + \maxavoid/2 + \auxw * sin(\auxa)});

        \coordinate (PD) at (1/2, {3/2*sqrt(3)});
        \coordinate (I5) at ($(PD) + (PA) - (I2)$);
        \coordinate (I6) at ($(PD) + (PA) - (I3)$);

        \coordinate (NN) at (1/2-\auxx, {sqrt(3)+\auxy});
        \coordinate (MM) at (1/2+\auxx, {sqrt(3)+\auxy});

        \draw[fill=black!5] (I5) -- (MM) -- (NN) -- (I6) -- (PD) -- cycle;

        \draw [{Bar[width=2mm]}-{Bar[width=2mm]}] ($(MM) + (0.2,0)$) -- ++(0,{sqrt(3)/2 - \auxy}) node[pos=0.5, right] {\tiny $\sqrt{27} /14$};
        
        \draw  (MM) -- (NN) node[pos=0.5, below] {\tiny  $1/7$};
        
    % \draw[line width=0.4pt] (MM) circle ({1/7});
    
        \draw  (I6) -- (NN) node[pos=0.5, left] {$s_2$};

        \draw  (I6) -- (PD) node[pos=0.4, above] {$s_3$};

        % \draw[very thin] (I5) circle[radius=0.5pt];
        % \draw[very thin] (MM) circle[radius=0.5pt];
\end{tikzpicture}
    \end{minipage}%
    \begin{minipage}[t]{0.35\textwidth}
    \centering
    \vspace{1.7em}
    \pgfmathsetmacro{\scale}{3}
\tikzset{every picture/.style={scale=\scale}}

\input{figures/secondcoloring/variables}

\newcommand{\circles}[1]{%
    \begin{scope}[shift={(#1)}]
    
    \draw[fill=black] (0,0) circle (0.01);
    \draw[dash dot, line width=0.1pt] (0,0) circle (\minavoid);
    \draw[dashed, line width=0.1pt] (0,0) circle (\maxavoid);
    \draw[dotted, line width=0.4pt] (0,0) circle (1);
        
    \end{scope}
}

\begin{tikzpicture}
        \coordinate (MAA) at (1/2, {- sqrt(3)/2 + \maxavoid});
        
        \coordinate (X) at (1/2, {sqrt(3)/2 - \maxavoid});
        
        \coordinate (XXX) at ({1/2 + sqrt(3)/2 - \maxavoid}, 0);
        \coordinate (XX) at ({1/2 - sqrt(3)/2 + \maxavoid}, 0);
        
        \draw[fill=black!5] (X)-- (XX) -- (MAA)  -- (XXX) -- cycle;

        \draw  (XX) -- (MAA) node[pos=0.5, below left] {\tiny $d_{\min}/\sqrt{2}$};
        
        \draw [{Bar[width=2mm]}-{Bar[width=2mm]}] ($(XX) + (0,0.3)$) -- ++({\minavoid},0) node[pos=0.5, above] {\tiny $d_{\min}$};
        
        \draw [{Bar[width=2mm]}-{Bar[width=2mm]}] ($(X) + (0.3,0)$) -- ++(0,{-\minavoid}) node[pos=0.5, right] {\tiny $d_{\min}$};

\end{tikzpicture}
    \end{minipage}%
    % \begin{minipage}[t]{0.6\textwidth}
    % \begin{align*}
    %     \alpha_1 & =  \\
    %     \alpha_2 & = \\
    %     \alpha_3 & = 
    %     % s_1 & = 1/7 \\
    %     % t_1 & = d_{\max} + \sqrt{3}/2 \\
    %     % s_2 & = \sqrt{ \frac{3 - d_{\max} \sqrt{3}}{2}
    %     % - \left(1 - \frac{d_{\max}}{2t} \right)
    %     %      \sqrt{
    %     %     \frac{1}{4} - \frac{4}{4t_1^2 + 1}
    %     %     - t_1^2 + \frac{15}{4}
    %     %     } } \\
    %     % s_3 & = ... 
    %     % h_1 & = \sqrt{27} /14
    % \end{align*}
    % \end{minipage}
    \caption{An axisymmetric pentagon and a square together are used for the red color avoiding points at distance $d$ in the second coloring. $s_2$ and $s_3$ are implicitly defined in \cref{fig:second_heptagon}.
    }
    \label{fig:second_pentagon}
\end{figure}

\begin{figure}[H]
    \centering
    \pgfmathsetmacro{\scale}{3}
\tikzset{every picture/.style={scale=\scale}}

\input{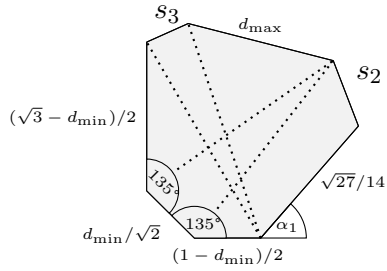}

\begin{tikzpicture}
        \coordinate (AA) at (1, 0);
        \coordinate (PA) at (1/2, {sqrt(3)/2});
        
        \coordinate (X) at (1/2, {sqrt(3)/2 - \maxavoid});
        
        \coordinate (XXX) at ({1/2 + sqrt(3)/2 - \maxavoid}, 0);
        
        \coordinate (W) at (3/4, {3*sqrt(3)/4 - \maxavoid/2});
        \coordinate (R) at (5/4, {sqrt(3)/4 + \maxavoid/2});
    
        \coordinate (I3) at ([shift=(180-\auxa:\auxw)]R);
        \coordinate (I4) at ([shift=(180-\auxa:-\auxw)]W);
        
        \coordinate (NN) at (1/2-\auxx, {sqrt(3)+\auxy});
        \coordinate (MM) at (1/2+\auxx, {sqrt(3)+\auxy});

        \pgfmathsetmacro{\yshift}{-sqrt(3)}
        \coordinate (MMM) at ([xshift=-1cm, yshift=\yshift cm)]MM);
        
        \pgfmathsetmacro{\yshift}{-sqrt(3)}
        \coordinate (NNN) at ([xshift=1cm, yshift=\yshift cm)]NN);
        
        \draw[fill=black!5] (X) -- (XXX) -- (AA) -- (NNN)  -- (I4) -- (I3) -- (PA) -- cycle;

        \draw[thick, dotted] (X) -- (I4);
        \draw[thick, dotted] (XXX) -- (I4);
        \draw[thick, dotted] (AA) -- (I3);
        \draw[thick, dotted] (AA) -- (PA);

        % \draw [{Bar[width=2mm]}-{Bar[width=2mm]}] ($(X) + (-0.1,0)$) -- ++(0,{-\minavoid/2}) node[pos=0.5, left] {\tiny $d_{\min}/2$};
        
        % \draw [{Bar[width=2mm]}-{Bar[width=2mm]}] ($(XXX) + (0,-0.1)$) -- ++({-\minavoid/2},0) node[pos=0.5, below] {\tiny $d_{\min}/2$};
        
        % \draw [{Bar[width=2mm]}-{Bar[width=2mm]}] ($(X) + (0,-0.5)$) -- ++({1/2},0) node[pos=0.5, below] {$1/2$};

        \draw  (XXX) -- (X) node[pos=0.5, below left] {\tiny $d_{\min}/\sqrt{2}$};

        \draw  (XXX) -- (AA) node[pos=0.5, below] {\tiny $(1 - d_{\min})/2$};

        \draw  (X) -- (PA) node[pos=0.5, left] {\tiny $(\sqrt{3} - d_{\min})/2$};

        \draw  (PA) -- (I3) node[pos=0.5, above] {$s_3$};
        
        \draw  (I4) -- (I3) node[pos=0.5, above] {\tiny $d_{\max}$};

        \draw  (I4) -- (NNN) node[pos=0.5, above right] {$s_2$};
        
        \draw  (AA) -- (NNN) node[pos=0.5, right] {\tiny $\sqrt{27}/14$};

        \draw[fill=black!5] (XXX) -- ++(0:0.14) arc (0:135:0.14) node[pos=0.55, below, inner sep=4pt] {\tiny $135\degree$}  -- cycle;
        
        \draw[fill=black!5] (X) -- ++(-45:0.14) arc (-45:90:0.14) node[pos=0.57, below, inner sep=3pt, rotate=-45] {\tiny $135\degree$}  -- cycle;
        
        \draw (AA) -- ++(0:0.2) arc (0:{45+acos(47/49)/4}:0.2) node[pos=0.5, below left, inner sep=1pt] {\tiny $\alpha_1$}  -- cycle;
        
        % \draw [{Bar[width=2mm]}-{Bar[width=2mm]}] ($(XXX) + (-0.7,0)$) -- ++(0,{sqrt(3)/2}) node[pos=0.5, left] {$\sqrt{3}/2$};
        
        % \draw [{Bar[width=2mm]}-{Bar[width=2mm]}] ($(AA) + (0.7,0)$) -- ++(0,{1.5}) node[pos=0.5, right] {$h_2$};
        
        % \draw [{Bar[width=2mm]}-{Bar[width=2mm]}] ($(X) + (0,0.95)$) -- ++(1.5,0) node[pos=0.5, above] {$w_1$};
\end{tikzpicture}
    \caption{A heptagon in which four of the diagonals have unit length, highlighted by dotted lines. Used for the orange, green, yellow, and turquoise color avoiding points at unit distance in the second coloring. We do not give a closed form solution for $s_2$ and $s_3$ but note that they are well defined. The angle $\alpha_1$ is defined in \cref{fig:second_hexagon}.}
    \label{fig:second_heptagon}
\end{figure}

\begin{figure}[H]
    \begin{minipage}[t]{0.5\textwidth}
    \centering
    \vspace{0.2em}
    \pgfmathsetmacro{\scale}{3}
\tikzset{every picture/.style={scale=\scale}}

\input{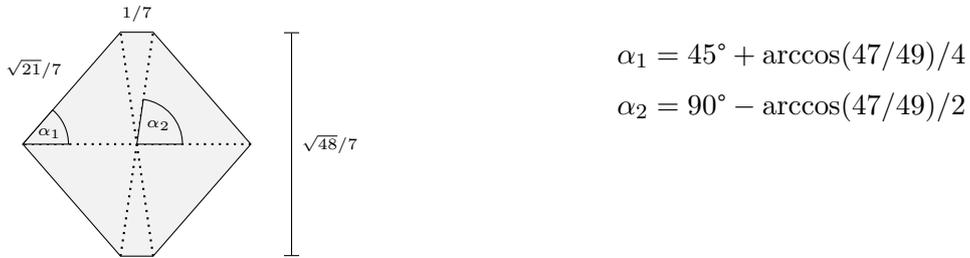}

\begin{tikzpicture}
        \coordinate (C) at (1, {sqrt(3)});
        \coordinate (B) at (0, {sqrt(3)});

        \coordinate (M) at (1/2-\auxx, {sqrt(3)-\auxy});
        \coordinate (N) at (1/2+\auxx, {sqrt(3)-\auxy});

        \coordinate (NN) at (1/2-\auxx, {sqrt(3)+\auxy});
        \coordinate (MM) at (1/2+\auxx, {sqrt(3)+\auxy});

        \draw[fill=black!5] (M) -- (N) -- (C) -- (MM) -- (NN) -- (B) -- cycle;

        \draw[thick, dotted] (C) -- (B);
        \draw[thick, dotted] (M) -- (MM);
        \draw[thick, dotted] (N) -- (NN);

        \node at ($(MM)!0.5!(NN)$) [above] {\tiny  $1/7$};

        \node at ($(B)!0.5!(NN)$) [above left] {\tiny $\sqrt{21}/7$};

        \draw [{Bar[width=2mm]}-{Bar[width=2mm]}] ($(M) + (0.75,0)$) -- ++(0,{2*\auxy}) node[pos=0.5, right] {\tiny  $\sqrt{48}/7$};

        \draw[fill=black!5] ({1/2},{sqrt(3)}) -- ++(0:0.2) arc (0:{90-acos(47/49)/2}:0.2) node[pos=0.45, below left, inner sep=1pt] {\tiny $\alpha_2$}  -- cycle;
        
        \draw[fill=black!5] (0,{sqrt(3)}) -- ++(0:0.2) arc (0:{45+acos(47/49)/4}:0.2) node[pos=0.5, below left, inner sep=1pt] {\tiny $\alpha_1$}  -- cycle;

        % \draw[dotted, line width=0.4pt] (B) circle ({sqrt(21)/7});
\end{tikzpicture}
    \end{minipage}%
    \begin{minipage}[t]{0.5\textwidth}
    \begin{align*}
        \alpha_1 & = 45\degree + \arccos(47/49)/4 \\
        \alpha_2 & = 90\degree - \arccos(47/49)/2
    \end{align*}
    \end{minipage}
    \caption{A centrosymmetric hexagon in which three of the diagonals have unit length, highlighted by dotted lines. Used for the blue color avoiding points at unit distance in the second coloring.
    }
    \label{fig:second_hexagon}
\end{figure}

\end{document}